\documentclass[12pt]{amsart}
\makeatletter
\let\over\@@over
\makeatother
\outer\def\proclaim #1. #2\par{\medbreak
  \noindent{\bfseries#1.\enspace}{\slshape#2\par}%
  \ifdim\lastskip<\medskipamount \removelastskip\penalty55\medskip\fi}
\let\i\subset

\begin{document}

\title{Littlewood-Richardson semigroups}
\author{Andrei Zelevinsky}
\address{\hskip-\parindent Andrei Zelevinsky,
Department of Mathematics,
Northeastern University,
Boston, MA 02115}

\email{andrei@neu.edu}

\thanks{This work is supported in part by NSF grant DMS-9625511;
research at MSRI is supported in part by NSF grant DMS-9022140.}

\begin{abstract}
This note is an extended abstract of my talk at 
the workshop on Representation Theory and Symmetric Functions,
MSRI, April 14, 1997.
We discuss the problem of finding an explicit description of 
the semigroup $LR_r$ of triples of partitions 
of length $\leq r$ such that the corresponding Littlewood-Richardson 
coefficient is non-zero. 
After discussing the history of the problem and previously known results, 
we suggest a new approach based on the ``polyhedral''
combinatorial expressions for the Littlewood-Richardson 
coefficients.
\end{abstract}

\maketitle

This note is an extended abstract of my talk at 
the workshop on Representation Theory and Symmetric Functions,
MSRI, April 14, 1997.
I thank the organizers (Sergey Fomin, Curtis Greene,
Phil Hanlon and Sheila Sundaram) for bringing together a group of 
outstanding combinatorialists and for giving me a chance to bring to their 
attention some of the problems that I find very exciting and beautiful.

\medskip

For $r \geq 1$, let 
$$P_r = \{\lambda = (\lambda_1, \ldots, \lambda_r) \in {\bf Z}^r:
\lambda_1 \geq \cdots \geq \lambda_r \geq 0\}$$
be the semigroup of partitions of length $\leq r$.
Our main object of study will be the set
$$LR_r =\{(\lambda,\mu,\nu): \lambda,\mu,\nu \in P_r, 
c^\lambda_{\mu \nu} > 0\} \ ,$$
where $c^\lambda_{\mu \nu}$ is the Littlewood-Richardson coefficient.
Recall that $P_r$ is the set of highest weights of polynomial irreducible
representations of $GL_r ({\bf C})$; if $V_\lambda$ is the 
irreducible representation of $GL_r ({\bf C})$ with highest weight $\lambda$
then $c^\lambda_{\mu \nu}$ is the multiplicity of $V_\lambda$ in 
$V_\mu \otimes V_\nu$.
Equivalently, the $c^\lambda_{\mu \nu}$ are the structure constants
of the algebra of symmetric polynomials in $r$ variables 
with respect to the basis of Schur polynomials.
We call $LR_r$ the {\it Littlewood-Richardson semigroup} of order $r$;
this name is justified by the following.

\proclaim Theorem 1. $LR_r$ is a finitely generated subsemigroup of 
the additive semigroup $P_r^3 \i {\bf Z}^{3r}$.

This is a special case of a much more general result 
well known to the experts in the invariant theory.
A short proof (valid for any reductive group instead of 
$GL_r ({\bf C})$) can be found in [E]; 
A.~Elashvili attributes this proof to M.~Brion and F.~Knop.
The semigroup property also follows at once from 
``polyhedral" expressions for $c^\lambda_{\mu \nu}$ 
that will be discussed later. 

\medskip

\noindent {\bf Problem A.} Describe $LR_r$ explicitly. 

\medskip

I have been interested in this problem for several years 
(e.g., in [BZ] the set $\{\lambda: (\lambda, \delta, \delta) \in LR_r\}$
was determined, where $\delta = (r-1, \ldots, 1, 0)$;
this proves a special case of Kostant's conjecture).
Practically nothing is known about the list of indecomposable generators
of $LR_r$ for general $r$.
We will discuss the ``dual'' approach, namely we would like to describe 
the facets of the polyhedral convex cone $LR_r^{\bf R} \i {\bf R}^{3r}$
generated by $LR_r$. 
A remarkable progress in this direction was recently made by 
A.~Klyachko in [K].
Before discussing his results, let us note that $c^\lambda_{\mu \nu}$
is given by the classical Littlewood-Richardson rule (see e.g., 
[M]), which in principle makes Problem~A purely combinatorial. 
In particular, the Littlewood-Richardson rule (or just the definition) 
readily implies the following properties of $LR_r$.

\medskip

\noindent {\bf Homogeneity.} $|\lambda| = |\mu| + |\nu|$ 
for $(\lambda, \mu, \nu) \in LR_r$, where 
$|\lambda| = \lambda_1 + \cdots + \lambda_r$.

\medskip

\noindent {\bf Stability.} $LR_{r+1} \cap {\bf Z}^{3r} = LR_r$,
where 
${\bf Z}^{3r} = \{(\lambda, \mu, \nu) \in {\bf Z}^{3(r+1)}: 
\lambda_{r+1} = \mu_{r+1} = \nu_{r+1} = 0\}$.
Even stronger, we have $LR_{r+1} \cap {\bf Z}^{3r+2} = LR_r$,
where 
${\bf Z}^{3r+2} = \{(\lambda, \mu, \nu) \in {\bf Z}^{3(r+1)}: 
\lambda_{r+1} = 0\}$.

\medskip

Littlewood-Richardson semigroups appear naturally in several
other contexts:

\smallskip

\noindent 1. Hall algebra, extensions of abelian $p$-groups: see [M].

\smallskip

\noindent 2. Schubert calculus on Grassmannians: see [F].

\smallskip

\noindent 3. Polynomial matrices and their invariant factors: see [T].

\smallskip

\noindent 4. Eigenvalues of sums of Hermitian matrices. 

\smallskip

Let us discuss the last item in more detail.
For a Hermitian matrix $A$ of order $r$, let $\lambda (A)$
denote the sequence of eigenvalues of $A$ arranged in a weakly 
decreasing order (recall that $A$ is Hermitian if $A^* = A$,
and such a matrix always has real eigenvalues). 
Let $HE_r$ denote the set of triples $(\lambda, \mu, \nu) \in {\bf R}^{3r}$
such that $\lambda = \lambda (A+B), \mu = \lambda (A)$, and 
$\nu = \lambda (B)$ for some Hermitian matrices $A$ and $B$ of order $r$. 
The following  counterpart of Theorem~1 for $HE_r$ is highly non-trivial.

\proclaim Theorem 2. $HE_r$ is a polyhedral convex cone in  ${\bf R}^{3r}$.

\medskip

\noindent {\bf Problem B.} Describe $HE_r$ explicitly. 

\medskip

Problems A and B are closely related to each other.
They have a long history.
Problem~B was probably first posed by I.M.~Gelfand in the late 40's
(eigenvalues of the sum of two Hermitian matrices
were studied already by H.~Weyl in 1912, but I believe 
I.M~Gelfand was the first who suggested to study the cone $HE_r$ 
as a whole rather than concentrate on individual eigenvalues).
A solution was announced by V.B.~Lidskii in [L1], but the details
of the proof were never published.
F.A.~Berezin and I.M.~Gelfand in [BG] discussed the relationships
between	Problems A and B; in particular, they suggested the following
remarkable equality:
$$HE_r \cap P_r^3 = LR_r. \eqno (1)$$
A.~Horn in [H] solved Problem~B for $r \leq 4$ and conjectured 
a general answer. 
To formulate his conjecture we need some terminology.
Let $[1,r]$ denote the set $\{1, 2, \ldots, r\}$.
For a subset $I = \{i_1 < i_2 < \cdots < i_s\} \i [1,r]$,
we denote by $\rho(I) \i P_s$ the partition 
$$\rho(I) = (i_s - s, \ldots, i_2 - 2, i_1 - 1).$$
A triple of subsets $I, J, K \i [1,r]$ will be called 
$HE$-{\it consistent} if they have the same cardinality $s$
and $(\rho(I), \rho(J), \rho(K)) \in HE_s$.
For $\lambda \in {\bf R}^r$ and $I \i [1,r]$, we will write
$|\lambda|_I = \sum_{i \in I} \lambda_i$;
in particular, 
$|\lambda|_{[1,r]} = |\lambda| = \lambda_1 + \cdots + \lambda_r$.

\medskip

\noindent {\bf Horn's Conjecture.} Let $\lambda, \mu$, and $\nu$ be 
vectors in ${\bf R}^r$ with weakly decreasing components.
Then $(\lambda,\mu,\nu) \in HE_r$ if and only if 
$|\lambda| = |\mu| + |\nu|$ and $|\lambda|_I \leq |\mu|_J + |\nu|_K$
for all $HE$-consistent triples $(I,J,K)$ in $[1,r]$. 

\medskip

The proofs of both Horn's Conjecture and (1) were announced by
B.V.~Lidskii (not to be confused with the author of [L1]) in [L2].
Unfortunately, as in the case of [L1], the detailed proofs of the results in
[L2] never appeared.
This justifies  A.~Klyachko's claim in [K] that even 
Theorem~2 has not been proved before.

\smallskip 

Let us now discuss the results in [K]. 
First of all, A.~Klyachko proves Theorem~2;
moreover, he gives the following description of the facets of $HE_r$,
which is very close (but not totally equivalent) to Horn's Conjecture.
Modifying the definition of $HE$-consistent triples.
we will call a triple of subsets $I, J, K \i [1,r]$ 
$LR$-{\it consistent} if they have the same cardinality $s$,
and $(\rho(I), \rho(J), \rho(K)) \in LR_s$.

\proclaim Theorem 3. Horn's conjecture becomes true if
$HE$-consistency in the formulation is replaced by $LR$-consistency.
Moreover, the inequalities $|\lambda|_I \leq |\mu|_J + |\nu|_K$
for all $LR$-consistent triples $(I,J,K)$ in $[1,r]$ are independent,
i.e., they correspond to facets of the polyhedral convex cone 
$HE_r$. 

A.~Klyachko also proves the following weaker version of (1).
Let $LR_r^{\bf Q}$ be the set of all linear combinations
of triples in $LR_r$ with positive rational coefficients;
equivalently, $LR_r^{\bf Q} = \cup_{N \geq 1}\,  {1 \over N} LR_r$.

\proclaim Theorem 4. $HE_r \cap {\bf Q}_{\geq 0}^{3r} = LR_r^{\bf Q}$.

Theorems 3 and 4 appear in [K] as a by-product 
of the study of stability criteria for toric vector bundles 
on the projective plane $P^2$. 
In view of these theorems, the equality (1) and Horn's Conjecture 
would follow from the affirmative answer to the following

\medskip

\noindent {\bf Saturation Problem.} Is it true that 
$LR_r^{\bf Q} \cap P_r^3 = LR_r$?

\medskip

In other words, does the fact that $c_{N\mu, N\nu}^{N\lambda} > 0$
for some $N \geq 1$ imply that $c_{\mu \nu}^\lambda > 0$?
This is true and easy to check for $r \leq 4$. 
On the other hand, an obvious analogue of the problem for type $B$
has negative answer (as pointed out to me by M.~Brion, counterexamples 
can be found in [E]). 

\medskip

\noindent {\bf Examples.} Here are the linear inequalities 
corresponding to $LR$-consistent triples for $r \leq 3$;
combined with the conditions $\lambda_1 \geq \cdots \geq \lambda_r,
\, \mu_1 \geq \cdots \geq \mu_r, \, \nu_1 \geq \cdots \geq \nu_r$,
and $|\lambda| = |\mu| + |\nu|$, they provide a description
of the cone $HE_r$. 

\begin{itemize} 

\item  $r = 1$: no inequalities;

\item  $r = 2$: $\lambda_1 \leq \mu_1 + \nu_1$, 
$\lambda_2 \leq {\rm min}\, (\mu_1 + \nu_2, \mu_2 + \nu_1)$;

\item  $r = 3$: $\lambda_1 \leq \mu_1 + \nu_1$, 
$\lambda_2 \leq {\rm min}\, (\mu_1 + \nu_2, \mu_2 + \nu_1)$,
$\lambda_3 \leq {\rm min}\, (\mu_1 + \nu_3, \mu_2 + \nu_2, \mu_3 + \nu_1)$,
$\lambda_1 + \lambda_2  \leq \mu_1 + \mu_2 + \nu_1 + \nu_2$, 
$\lambda_1 + \lambda_3 \leq {\rm min}\, (\mu_1 + \mu_2 + \nu_1 + \nu_3,
\mu_1 + \mu_3 + \nu_1 + \nu_2)$,
$\lambda_2 + \lambda_3 \leq {\rm min}\, (\mu_1 + \mu_2 + \nu_2 + \nu_3,
\mu_1 + \mu_3 + \nu_1 + \nu_3, \mu_2 + \mu_3 + \nu_1 + \nu_2)$.
For instance, the inequality 
$\lambda_2 + \lambda_3 \leq \mu_1 + \mu_3 + \nu_1 + \nu_3$ 
corresponds to the triple $(I,J,K) = (\{2,3\}, \{1,3\}, \{1,3\})$,
which is $LR$-consistent because the triple of partitions 
$(\rho(I), \rho(J), \rho(K)) = ((1,1), (1,0), (1,0))$ obviously 
belongs to $LR_2$. 
\end{itemize}

Assuming the affirmative answer in the Saturation Problem,
Theorem~3 provides a recursive procedure for describing the semigroup $LR_r$. 
Although quite elegant, this procedure is not very explicit from combinatorial
point of view. 
Thus, we would like to formulate the following

\medskip

\noindent {\bf Problem C.} Find a non-recursive description of $LR_r$.

\medskip

Equivalently, Problem~C asks for a non-recursive description of 
$LR$-consistent triples. 
We would like to suggest an elementary combinatorial approach to this problem
based on the ``polyhedral'' expressions for the coefficients 
$c^\lambda_{\mu \nu}$ given in [BZ]. 
To present such an  expression, it will be convenient to modify 
Littlewood-Richardson coefficients as follows. 
We will consider triples $(\bar \lambda, \bar \mu, \bar \nu)$
of dominant integral weights for the group $SL_r$. 
Let $V_{\bar \lambda}$ be the irreducible $SL_r$-module with highest weight
$\bar \lambda$, and let $c_{\bar \lambda \bar \mu \bar \nu}$
denote the dimension of the space of $SL_r$-invariants in the triple 
tensor product $V_{\bar \lambda} \otimes V_{\bar \mu} \otimes V_{\bar \nu}$.
The relationship between the $c_{\bar \lambda \bar \mu \bar \nu}$
and the Littlewood-Richardson coefficients is as follows. 
We will write each of the weights $\bar \lambda, \bar \mu$ and  $\bar \nu$
as a nonnegative integer linear combination of fundamental weights
$\omega_1, \omega_2, \ldots, \omega_{r-1}$ (in the standard numeration):
$$\bar \lambda = l_1 \omega_1 + \cdots + l_{r-1} \omega_{r-1}, \,\,
\bar \mu = m_1 \omega_1 + \cdots + m_{r-1} \omega_{r-1}, \,\,
\bar \nu = n_1 \omega_1 + \cdots + n_{r-1} \omega_{r-1}. \eqno (2)$$
The definitions readily imply that if $\lambda, \mu, \nu \in P_r$ are such that
$|\lambda| = |\mu| + |\nu|$ then 
$c^\lambda_{\mu \nu} = c_{\bar \lambda \bar \mu \bar \nu}$,
where the coordinates $l_s, m_s$ and $n_s$ in (2) are given by
$$l_s = \lambda_{r-s} - \lambda_{r-s+1}, \,\, m_s = \mu_s - \mu_{s+1}, \,\,
n_s = \nu_s - \nu_{s+1}. \eqno (3)$$
Thus, the knowledge of $LR_r$ is equivalent to the knowledge
of the semigroup 
$$\overline {LR}_r =  \{(\bar \lambda, \bar \mu, \bar \nu) 
\in {\bf Z}_{\geq 0}^{3(r-1)}: c_{\bar \lambda \bar \mu \bar \nu} > 0\}.$$

Passing from $LR_r$ to $\overline {LR}_r$ has two important advantages.
First, the coefficients $c_{\bar \lambda \bar \mu \bar \nu}$ are more 
symmetric than the original Littlewood-Richardson coefficients:
they are invariant under the $12$-element group generated by all 
permutations of three weights $\bar \lambda, \bar \mu$ and  $\bar \nu$,
together with the transformation replacing each of these weights 
with its dual (i.e., sending $(l_s, m_s, n_s)$ to 
$(l_{r-s}, m_{r-s}, n_{r-s})$). 
Second, the dimension of the ambient space reduces by 2, from 
$3r-1$ to $3(r-1)$. 
On the other hand, $\overline {LR}_r$ has at least one potential disadvantage:
the condition $|\lambda| = |\mu| + |\nu|$ is replaced by a more complicated
condition that $\sum_s s(l_s + m_s + n_s)$ is divisible by $r$
(in more invariant terms, this means that $\bar \lambda + \bar \mu + \bar \nu$
must be a radical weight, i.e., belongs to the root lattice). 
To illustrate both phenomena, one can compare the description of $LR_2$
given above with the following description of $\overline {LR}_2$ 
which is equivalent to the classical Clebsch-Gordan rule:
$\overline {LR}_2$ consists of triples of nonnegative integers 
$(l_1, m_1, n_1)$ satisfying the triangle inequality and such that 
$l_1 + m_1 + n_1$ is even.

Let us now give a combinatorial expression for
$c_{\bar \lambda \bar \mu \bar \nu}$ (this is one of several 
such expressions found in [BZ]).
Consider a triangle in ${\bf R}^2$, and subdivide it into small 
triangles by dividing each side into $r$ equal parts and joining 
the points of the subdivison by the line segments parallel to the sides
of our triangle. 
Let $Y_r$ denote the set of all vertices of the small triangles,
with the exception of the three vertices of the original triangle.
Introducing barycentric coordinates, we identify
$Y_r$ with the set of integer triples $(i,j,k)$ such that 
$0 \leq i, j, k < r$ and $i+j+k = r$.
Let ${\bf Z}^{Y_r}$ be the set of integer families
$(y_{ijk})$ indexed by $Y_r$; we think of $y \in {\bf Z}^{Y_r}$
as an integer ``matrix'' with $Y_r$ as the set of ``matrix positions.''
To every $y \in {\bf Z}^{Y_r}$ we associate the partial line sums
$$l_{ts} (y) = \sum_{j=t}^s y_{r-s,j,s-j}, \,\, 
m_{ts} (y) = \sum_{k=t}^s y_{s-k,r-s,k}, \,\,
n_{ts} (y) = \sum_{i=t}^s y_{i,s-i,r-s}, \eqno (4)$$
where $0 \leq t \leq s \leq r$.
We call these linear forms on ${\bf R}^{Y_r}$ {\it tails}, and we say that 
$y \in {\bf R}^{Y_r}$ is {\it tail-positive} if all tails of $y$ are $\geq 0$.
We also say that a linear form on ${\bf R}^{Y_r}$ is {\it tail-positive}
if it is a nonnegative linear combination of tails. 

\proclaim Theorem 5. For any triple $(\bar \lambda, \bar \mu, \bar \nu)$
as in (2), the coefficient $c_{\bar \lambda \bar \mu \bar \nu}$
is equal to the number of tail-positive $y \in {\bf Z}^{Y_r}$
with prescribed values of line sums
$$l_{0s} (y) = l_s, \,\, m_{0s} (y) = m_s, \,\, n_{0s} (y) = n_s
\quad \quad (1 \leq s \leq r-1). \eqno (5)$$

In other words, let $T_r \i {\bf Z}^{Y_r}$ denote the semigroup of 
tail-positive elements, and let 
$\sigma:{\bf Z}^{Y_r} \to {\bf Z}^{3(r-1)}$ be the projection given by (5).
Then Theorem~5 says that
$$\sigma (T_r) =   \overline {LR}_r. \eqno (6)$$
In particular, this implies at once that $\overline {LR}_r$
(and hence $LR_r$) is a semigroup.
Furthermore, Theorem~5 implies the following description 
of the convex cone $\overline {LR}^{\bf R}_r$ generated 
by $\overline {LR}_r$. 

\proclaim Corollary 6. A linear form $f$ on ${\bf R}^{3(r-1)}$
takes nonnegative values on $\overline {LR}^{\bf R}_r$ if and only if
the form $f \circ \sigma$ on ${\bf R}^{Y_r}$ is tail-positive.

Returning to the Littlewood-Richardson semigroup $LR_r$, 
we have the projection $\partial:LR_r \to \overline {LR}_r$
given by (3). 
This projection extends by linearity to a projection
$\partial: {\bf R}^{3r-1} \to {\bf R}^{3(r-1)}$, where 
${\bf R}^{3r-1}$ is the subspace of triples 
$(\lambda, \mu, \nu) \in {\bf R}^{3r}$ satisfying 
$|\lambda| = |\mu| + |\nu|$. 
It is clear that the cone ${LR}^{\bf R}_r \i {\bf R}^{3r-1}$ 
is given by the linear inequalities $f \circ \partial \geq 0$
for all linear forms $f$ as in Corollary~6. 
This suggests the following strategy for determining the set of 
$LR$-consistent triples. 
Take a triple of subsets $(I,J,K)$ in $[1,r]$ of the same cardinality $s$,
consider the corresponding linear form 
$|\mu|_J + |\nu|_K - |\lambda|_I$ on ${\bf R}^{3r-1}$,
write this form as $f \circ \partial$, and compute the form 
$f \circ \sigma$ on ${\bf R}^{Y_r}$. 
A straightforward calculation gives
$$(f \circ \sigma)(y) = 
\sum_{(i,j,k) \in Y_r} \bigl(\# (I_{> i}) - \# (J_{> r-j}) 
- \# (K_{> r-k})\bigr) \,  y_{ijk}, \eqno (7)$$
where $\# (I_{> i})$ stands for the number of elements of $I$ which are 
$> i$.  
Taking into account Theorem~3, we obtain the following criterion
for $LR$-consistency.

\proclaim Theorem 7. A triple of subsets $(I,J,K)$ of the same cardinality
$s$ in $[1,r]$ is $LR$-consistent if and only if 
$|\rho(I)| = |\rho(J)| +|\rho(K)|$ and the form in (7) is tail-positive.

In particular, since every tail-positive linear form is obviously a 
nonnegative linear combination of the $y_{ijk}$, we obtain the following 
necessary condition for $LR$-consistency.

\proclaim Corollary 8. If a triple of subsets $(I,J,K)$  in $[1,r]$ is 
$LR$-consistent then 
$$\# (I_{> i}) \geq \# (J_{> r-j}) + \# (K_{> r-k}) \eqno (8)$$
for all $(i,j,k) \in Y_r$. 

It would be interesting to deduce this corollary directly from the 
Littlewood-Richardson rule.
One can show that (8) is not sufficient for $LR$-consistency.
In fact, Theorem~7 can be used to produce other necessary conditions for 
$LR$-consistency.
One can hope to solve Problem~C by generating  a system of necessary and 
sufficient conditions for $LR$-consistency using this method.

\end{document}